\documentclass[11pt,letterpaper]{article}
\usepackage[T1]{fontenc}
\usepackage{lmodern,amsmath,amssymb,amsthm,microtype}
\usepackage[margin=1in]{geometry}
\usepackage[hidelinks,bookmarksnumbered=true]{hyperref}
\usepackage{enumitem}
\numberwithin{equation}{section}
\newtheorem{theorem}{Theorem}[section]
\newtheorem{proposition}[theorem]{Proposition}
\newtheorem{lemma}[theorem]{Lemma}

\newcommand{\R}{\mathbb R}

\newcommand{\ii}{\mathrm i}
\newcommand{\dd}{\,\mathrm d}
\newcommand{\wt}{\widehat}
\newcommand{\cH}{\mathcal H}
\newcommand{\cQ}{\mathcal Q}
\newcommand{\re}{\operatorname{Re}}
\newcommand{\e}{\mathrm e}
\setlist[enumerate]{topsep=3pt,itemsep=3pt}
\title{Arithmetic truncations of a theta correlation:\\
Fourier positivity and localization of critical-line zeros}
\author{Wei Sun\\[3pt]
  \small Department of Mathematics and Statistics\\
  \small Concordia University\\
  \small Montreal, H3G 1M8, Canada\\
  \small\href{mailto:wei.sun@concordia.ca}{\texttt{wei.sun@concordia.ca}}}
\date{}
\makeatletter
\renewcommand{\@maketitle}{%
  \begin{center}
    {\LARGE\@title\par}
    \vspace{0.8em}
    {\large\begin{tabular}{c}\@author\end{tabular}\par}
  \end{center}}
\makeatother
\hypersetup{pdftitle={Arithmetic truncations of a theta correlation: Fourier positivity and localization of critical-line zeros},pdfauthor={Wei Sun}}

\begin{document}
\maketitle
\begin{abstract}
We truncate a weighted correlation of the Riemann theta kernel by retaining
the pairs of theta indices with product at most \(N\). For each tilt \(0<\tau\le\tfrac12\),
every mode of the tilt derivative with index \(k\ge2\) has a strictly
positive Fourier transform. For each fixed \(N\ge1\), the transform of
the sum over \(1\le k\le N\) is negative for all sufficiently large
\(\lvert t\rvert\). At tilt
\(\tfrac12\) we locate the end of its first interval of Fourier nonnegativity
at \(4N+(2/\pi)\log N+O(\log\log N)\). Assuming the Riemann hypothesis, the
same two leading terms hold with error \(o(\log N)\), uniformly when the
tilt ranges over a compact subinterval of \((0,\tfrac12]\). At zero tilt,
the Fourier transform of each finite truncation is negative at every
critical-line zero. For a fixed such zero of multiplicity \(m\), we locate
both endpoints of the surrounding negative interval and show that its width
has logarithm \(-\pi N/m+(5/(4m))\log N+o(\log N)\). This result is
unconditional. The leading term for the width depends on the divisor weights
of the omitted modes and the first nonzero derivative of \(\xi\) at the zero.
\end{abstract}

\noindent\textit{Research collaboration.} This paper is based on research
that Wei Sun and ChatGPT (Alex), an AI system developed by OpenAI, conducted
jointly from 17 August through 25 November 2026. Their respective
contributions are described at the end of the paper.

\smallskip
\noindent\textit{Keywords.} Riemann xi function, theta correlation, arithmetic
truncation, Fourier positivity, modified Bessel function, zero multiplicity.

\smallskip
\noindent\textit{2020 Mathematics Subject Classification.}
Primary 11M26. Secondary 42A38, 33C10.

\section{Introduction and main results}\label{sec:intro}

The Riemann xi function is the entire continuation of
\(\xi(s)=\tfrac12s(s-1)\pi^{-s/2}\Gamma(s/2)\zeta(s)\) from \(\re s>1\).
It satisfies \(\xi(0)=\xi(1)=\tfrac12\). The Riemann hypothesis (RH)
asserts that its zeros lie on \(\re s=\tfrac12\).

For \(Y>0\) and \(x\in\R\), define the theta sum \(\psi\) and its
associated kernel \(\phi\) by
\begin{equation}\label{eq:phi}
 \psi(Y)=\sum_{n\ge1}\e^{-\pi n^2Y},\qquad
 g(x)=\e^{x/4}\psi(\e^x),\qquad
 \phi(x)=4g''(x)-\tfrac14g(x).
\end{equation}
For \(0\le\tau\le\tfrac12\), define the weighted correlation \(H_\tau\)
and its tilt derivative \(Q_\tau\) by
\begin{equation}\label{eq:HQ}
 H_\tau(x)=\int_{\R}\cosh(\tau y)\phi(x+y)\phi(x-y)\dd y,\qquad
 Q_\tau(x)=\partial_\tau H_\tau(x).
\end{equation}
Write \(\wt f(t)=\int_{\R}f(x)\e^{-\ii tx}\dd x\). At
\(s=\tfrac12+\tau+\ii t\), \(t\in\R\), the two Fourier transforms satisfy
\begin{equation}\label{eq:completion}
 \wt H_\tau(t)=2|\xi(s)|^2,\qquad
 W_\tau(t):=\wt Q_\tau(t)=4\re\bigl(\xi'(s)\overline{\xi(s)}\bigr).
\end{equation}
By symmetry in \(y\), we may replace \(\cosh(\tau y)\) by \(\e^{\tau y}\).
With \(u=x+y,\ v=x-y\), the two Mellin factors are \(2\xi(\bar s)\)
and \(2\xi(1-s)\) (Subsection~\ref{sec:modes}). Since
\(\dd x\dd y=\tfrac12\dd u\dd v\), conjugation and the functional equation
give \(\wt H_\tau(t)=2|\xi(s)|^2\). The decay of \(\phi\) permits
differentiation in \(\tau\) at fixed \(t\), proving the second identity.

Our motivation comes from Sun's work on the modulus of \(\xi\). A
representation of \(|\xi|\) gave an analytic inequality equivalent to
RH~\cite{SunInequality}. A subsequent representation in terms of the Jacobi
theta function led to a family of polynomial inequalities whose joint
validity is equivalent to RH~\cite{SunModulus}. Here we study how the
signs of the transforms in~\eqref{eq:completion} change when the double
theta series for \(H_\tau\) is truncated.

\subsection{Arithmetic truncation and sign intervals}\label{sec:truncation}
Applying the differential operator in~\eqref{eq:phi} to each term of
\(g\) gives \(\phi(u)=\sum_{n\ge1}\phi_n(u)\), where
\[
 \phi_n(u)=2\pi n^2\e^{5u/4}(2\pi n^2\e^u-3)\e^{-\pi n^2\e^u},
 \qquad u\in\R .
\]
For each integer \(k\ge1\) and \(0\le\tau\le\tfrac12\), define the
\(k\)-th arithmetic mode of \(H_\tau\) on \(x\ge0\) by
\[
 H_{\tau,k}(x):=\int_\R\cosh(\tau y)
       \sum_{n\mid k}\phi_n(x+y)\phi_{k/n}(x-y)\dd y .
\]
For each positive divisor \(n\) of \(k\), set \(m=k/n\). The sum then
contains each ordered pair \((n,m)\) of positive integers with \(nm=k\)
exactly once. Interchanging \(n\) and \(m\) shows that the divisor sum
is even in \(y\). For \(x<0\), set \(H_{\tau,k}(x)=H_{\tau,k}(-x)\).
Grouping the terms of~\eqref{eq:HQ} on \(x\ge0\) and using the evenness
of \(H_\tau\) gives \(H_\tau=\sum_{k\ge1}H_{\tau,k}\) on \(\R\).
With \(Q_{\tau,k}=\partial_\tau H_{\tau,k}\), retaining the modes with
\(k=nm\le N\) gives the finite transforms
\begin{equation}\label{eq:finite}
 \cH_{\tau,N}(t)=\sum_{k=1}^N\wt H_{\tau,k}(t),\qquad
 \cQ_{\tau,N}(t)=\sum_{k=1}^N\wt Q_{\tau,k}(t),\qquad N\ge1 .
\end{equation}
All logarithms are natural. For \(G\ge0\), the notation \(F\ll_\epsilon G\)
means that \(|F|\le C_\epsilon G\), while a lower bound \(F\gg_a G\)
means that \(F\ge c_aG\). The constants are positive and uniform in
the stated ranges. Subscripts indicate the parameters on which the
constants may depend.

For positive tilt, how far from frequency zero does \(\cQ_{\tau,N}\)
remain nonnegative? For \(0<\tau\le\tfrac12\), let
\[
 \rho_N(\tau)=\sup\bigl(\{0\}\cup
 \{T>0:\ \cQ_{\tau,N}(t)\ge0\text{ for all }0\le t\le T\}\bigr).
\]
\begin{theorem}[First interval of Fourier nonnegativity]\label{th:window}
The numbers \(\rho_N(\tau)\) are finite and nondecreasing in \(N\).
\begin{enumerate}[label=\textup{(\roman*)}]
\item Unconditionally, as \(N\to\infty\),
\begin{equation}\label{eq:endpoint}
 \rho_N(1/2)=4N+\frac{2}{\pi}\log N+O(\log\log N).
\end{equation}
\item Assuming RH, for every fixed \(0<a<\tfrac12\),
\begin{equation}\label{eq:fixedstrip}
 \lim_{N\to\infty}\ \sup_{a\le\tau\le1/2}
 \frac{\big|\rho_N(\tau)-4N-(2/\pi)\log N\big|}{\log N}=0 .
\end{equation}
\end{enumerate}
\end{theorem}

Let \(N_{\min}(T)\) be the least cutoff for which \(\cQ_{1/2,N}\) is
nonnegative on all of \([0,T]\). Inverting the unconditional
endpoint estimate~\eqref{eq:endpoint} gives
\begin{equation}\label{eq:inverse}
 N_{\min}(T)=\frac T4-\frac{1}{2\pi}\log T+O(\log\log T)
 \qquad(T\to\infty).
\end{equation}
A negative value at a single test frequency gives an upper bound for
\(\rho_N\). The lower bound in Theorem~\ref{th:window} requires control
of the omitted modes at every earlier frequency.

At tilt zero, \(Q_0\) vanishes identically by symmetry in \(\tau\).
We therefore turn to the correlation itself. If \(\tfrac12+\ii\gamma\)
is a zero of \(\xi\), then \(\wt H_0(\gamma)=0\) by~\eqref{eq:completion}.
How wide is the interval about \(\gamma\) on which the finite transform
\(\cH_{0,N}\) is negative? Write \(d(k)=\sum_{j\mid k}1\) for the
divisor function and put
\begin{equation}\label{eq:AN}
 A_N=16\pi^3\sum_{k>N}d(k)k^{5/2}\e^{-2\pi k}.
\end{equation}

\begin{theorem}[The negative interval at a fixed zero]\label{th:zero}
Let \(s_\gamma=\tfrac12+\ii\gamma\) be a fixed zero of \(\xi\) of multiplicity
\(m\), and set \(C_\gamma=2|\xi^{(m)}(s_\gamma)/m!|^2\). For every \(N\ge1\),
\(\cH_{0,N}(\gamma)<0\). For all sufficiently large \(N\), the connected
component of \(\{t\in\R:\cH_{0,N}(t)<0\}\) containing \(\gamma\) is a bounded
interval \(I_N(\gamma)=(\gamma-h_N^-,\,\gamma+h_N^+)\), with
\begin{equation}\label{eq:width}
 h_N^- \sim (A_N/C_\gamma)^{1/(2m)},\qquad
 h_N^+ \sim (A_N/C_\gamma)^{1/(2m)}.
\end{equation}
In particular, the length \(|I_N(\gamma)|=h_N^-+h_N^+\) of this interval satisfies
\begin{equation}\label{eq:widthlog}
 \log |I_N(\gamma)|=-\frac{\pi}{m}N+\frac{5}{4m}\log N
                         +o_\gamma(\log N).
\end{equation}
These assertions require only the indicated zero, with no assumption on the
locations or multiplicities of other zeros.
\end{theorem}

For example, the logarithm of the width divided by \(N\) tends to
\(-\pi\) at a simple zero and to \(-\pi/2\) at a double zero if one exists.
The factor \(A_N\) retains the arithmetic variation just beyond the cutoff.
For two fixed simple zeros, this common factor cancels in~\eqref{eq:width},
giving the derivative comparison
\[
 \frac{|I_N(\gamma)|}{|I_N(\gamma')|}\longrightarrow
       \frac{|\xi'(s_{\gamma'})|}{|\xi'(s_\gamma)|} .
\]

\subsection{The completed transform and its arithmetic modes}\label{sec:modes}

For \(x\ge0\), set \(Y=\e^x\). Each summand \(\phi_n(x)\) in
Subsection~\ref{sec:truncation} is positive, since \(2\pi n^2Y>3\).
For each fixed integer \(j\ge0\), summing the differentiated terms gives
\(\lvert\phi^{(j)}(x)\rvert\le P_j(Y)\e^{-\pi Y}\) for some polynomial
\(P_j\), since \(\e^{-\pi(n^2-1)Y}\le\e^{-\pi(n^2-1)}\) for \(Y\ge1\).
Poisson summation gives \(1+2\psi(Y)=Y^{-1/2}(1+2\psi(1/Y))\) for \(Y>0\), hence
\(g(-x)-g(x)=\sinh(x/4)\). The operator in~\eqref{eq:phi} annihilates
this difference, so \(\phi\) is even and its derivatives have the same
double-exponential decay as \(x\to-\infty\).

For \(\re s>1\), integration by parts twice gives
\begin{align}
 \int_\R\phi(x)\e^{(s-1/2)x/2}\dd x
 &=s(s-1)\int_0^\infty\psi(Y)Y^{s/2}\,\frac{\dd Y}{Y}\notag\\
 &=s(s-1)\pi^{-s/2}\Gamma(s/2)\zeta(s)=2\xi(s).\label{eq:mellin}
\end{align}
The theta identity gives \(g,g'=O(\e^{-x/4})\) as \(x\to-\infty\), so the
boundary terms vanish. The decay of \(\phi\) makes its Mellin integral
entire in \(s\) and justifies Fubini's theorem and differentiation in
\(\tau\) for the correlation. Since \(\phi\) is real and even,
\(\xi(\bar s)=\overline{\xi(s)}\) and \(\xi(1-s)=\xi(s)\).
For \(s=\tfrac12+\tau+\ii t\), the substitution following~\eqref{eq:completion}
gives \(\tau y-\ii tx=(\bar s-\tfrac12)u/2+(\tfrac12-s)v/2\).
Hence~\eqref{eq:mellin} gives \(\wt H_\tau(t)=2\xi(\bar s)\xi(1-s)=2|\xi(s)|^2\),
including at zeros of \(\xi\).

For \(z>0\), let \(K_\nu\) be the modified Bessel function of the second
kind and define \(B_\nu\) by
\begin{equation}\label{eq:B}
 K_\nu(z)=\int_0^\infty \e^{-z\cosh u}\cosh(\nu u)\dd u,\qquad
 B_\nu(z)=(z^2+9)K_\nu(z)+6zK_\nu'(z).
\end{equation}
A prime on \(K_\nu\) denotes an argument derivative. For \(x\ge0\) set
\begin{equation}\label{eq:ak}
 a_k(\nu)=k^2\sum_{j\mid k}(k/j^2)^\nu,\qquad
 G_{\nu,k}(x)=\e^{5x/2}B_\nu(2\pi k\e^x),
\end{equation}
extending \(G_{\nu,k}\) evenly. We claim that
\begin{equation}\label{eq:mode}
 H_{\tau,k}=8\pi^2a_k(\tau)G_{\tau,k},\qquad
 Q_{\tau,k}=8\pi^2\{a_k'(\tau)G_{\tau,k}
                           +a_k(\tau)\partial_\tau G_{\tau,k}\}.
\end{equation}
For fixed \(k\), the divisor sum defining \(H_{\tau,k}\) is even in \(y\)
by the interchange \(n\leftrightarrow m\). We may therefore replace
\(\cosh(\tau y)\) by \(\e^{\tau y}\) before integrating the individual
terms. For the term \(\phi_n(x+y)\phi_m(x-y)\), where \(nm=k\), put
\(Y=\e^x,\ A=\pi n^2Y,\ C=\pi m^2Y,\ z=2\sqrt{AC}=2\pi kY\). The change
\(u=y+\tfrac12\log(A/C)\) gives
\[
 A\e^y+C\e^{-y}=z\cosh u,\qquad
 (2A\e^y-3)(2C\e^{-y}-3)=z^2+9-6z\cosh u.
\]
Using \(\e^{\tau y}=(m/n)^\tau\e^{\tau u}\) and integrating over \(\R\)
yields \(8\pi^2Y^{5/2}(mn)^2(m/n)^\tau B_\tau(z)\).
Summing over \(n\mid k\), with \(m=k/n\), gives \(H_{\tau,k}\)
in~\eqref{eq:mode}, and its \(\tau\)-derivative gives \(Q_{\tau,k}\).

To justify the grouping, after this change of variables the absolute
integrand is bounded by a polynomial in \(k,Y,\cosh u,|u|+\log(2k)\) times
\(\e^{|u|/2-z\cosh u}\). Fixed orders of differentiation in \(x,\tau\)
only change the polynomial. The estimate \(\cosh u\ge1+u^2/2\) bounds its
integral by a polynomial in \(k,Y\) times \(\e^{-2\pi kY}\).
The number of pairs is at most \(k\). Hence the half-line series and their
fixed derivatives converge absolutely and locally uniformly, up to \(x=0+\),
and in \(L^1(0,\infty)\). Indeed the integral of any fixed power of \(Y\)
against \(\e^{-2\pi kY}\dd Y/Y\) on \(Y\ge1\) decays exponentially in \(k\).
These facts justify termwise transformation and right differentiation.

The correlation in~\eqref{eq:HQ} also connects with earlier work on
Fourier kernels. Csordas \cite[\S3]{Csordas} obtains moment kernels
from the squared modulus of a Fourier transform by passing to sum and
difference variables. Expanding \(\cosh(\tau y)\) in~\eqref{eq:HQ} gives
these moment kernels as the coefficients of \(\tau^{2n}/(2n)!\).
Haglund \cite{Haglund}
studies the zeros of finite sums in an incomplete-gamma expansion of
\(\xi\). Our truncation acts on the two theta indices in the correlation
for \(2|\xi|^2\). Its tilt derivative connects with the horizontal
monotonicity of \(|\xi|\) studied by Sondow and Dumitrescu \cite{SD}
and by Matiyasevich, Saidak and Zvengrowski \cite{MSZ}.
The weight \(a_k(\nu)/k^2\) in~\eqref{eq:ak} is the generalized divisor
sum that accompanies \(K_\nu\) in the Fourier expansion of the real
analytic Eisenstein series \cite[(38)--(39)]{Zagier}.

The remainder of this paper is organized as follows.
Section~\ref{sec:shape} proves Fourier positivity for the individual
modes \(k\ge2\) and eventual negativity for each finite sum of tilt
derivatives. Section~\ref{sec:tails} estimates the omitted modes.
Section~\ref{sec:local} uses the untilted tail to determine the two ends
of the negative interval at a fixed zero, proving Theorem~\ref{th:zero}.
Section~\ref{sec:complete} bounds the complete tilt derivative at the
endpoint unconditionally and on compact positive-tilt strips under RH.
Section~\ref{sec:window} compares these bounds with the omitted tail to
prove Theorem~\ref{th:window} and the inverse cutoff formula~\eqref{eq:inverse}.
The final two unnumbered sections describe the formal verification and
the research contributions.

\section{Convexity and the exceptional first mode}\label{sec:shape}

We first give estimates for the integral in~\eqref{eq:B}. For \(z\ge1\) and
\(0\le\nu\le\tfrac12\),
\begin{equation}\label{eq:Kbound}
 |K_\nu(z)|+|K_\nu'(z)|+|\partial_\nu K_\nu(z)|
             +|\partial_\nu K_\nu'(z)|\ll z^{-1/2}\e^{-z}.
\end{equation}
In addition,
\begin{equation}\label{eq:Kasym}
 \e^z\sqrt z\,K_\nu(z)\longrightarrow\sqrt{\pi/2},\qquad
 -\e^z\sqrt z\,K_\nu'(z)\longrightarrow\sqrt{\pi/2}
\end{equation}
uniformly in this range of orders. To see both assertions, put
\(u=v/\sqrt z\) in the integral. The exponential after removal of
\(\e^{-z}\) is bounded by \(\e^{-v^2/2}\). A derivative in \(z\) introduces
\(-\cosh u\) and an order derivative introduces \(u\sinh(\nu u)\) in
place of \(\cosh(\nu u)\). The resulting integrands have an integrable
majorant \(C(1+v)\e^{2v-v^2/2}\). For~\eqref{eq:Kasym} they converge
uniformly in \(\nu\) at each \(v\) to \(\e^{-v^2/2}\), after use of the
indicated sign. This proves both statements and the leading Bessel
asymptotic used below. See also \cite[\S\S10.32, 10.40]{DLMF}.

Set \(D=z\partial_z\) and \(L=D+5/2\). Two integrations by parts in \(u\)
give \(D^2 K_\nu=(z^2+\nu^2)K_\nu\). Indeed
\((D^2-z^2)\e^{-z\cosh u}=\partial_u^2\e^{-z\cosh u}\), and both first
derivatives in the boundary terms vanish at zero. The variable \(x\) in
\eqref{eq:ak} gives \(\partial_x^jG_{\nu,k}(x)=\e^{5x/2}(L^jB_\nu)(2\pi k\e^x)\)
on the right half-line.

\begin{lemma}\label{lem:conv}
For \(z\ge4\pi\) and \(0\le\nu\le\tfrac12\), one has \(L^2B_\nu(z)>0\).
For \(0<\nu\le\tfrac12\), one also has \(\partial_\nu L^2B_\nu(z)>0\).
\end{lemma}
\begin{proof}
Write \(r=-DK_\nu/K_\nu\) and \(p=\partial_\nu r\). The integral shows
\(r>z\), since \(r/z\) is the average of \(\cosh u\) for a strictly
positive density on \(u>0\). The differential equation gives
\(Dr=r^2-z^2-\nu^2\). Consequently \(D(pK_\nu^2)=-2\nu K_\nu^2\), so
\[
 p(z)=\frac{2\nu}{K_\nu(z)^2}\int_z^\infty K_\nu(v)^2\,\frac{\dd v}{v}.
\]
The boundary term vanishes by~\eqref{eq:Kbound}, applied to
\(pK_\nu^2=-(D\partial_\nu K_\nu)K_\nu+(DK_\nu)\partial_\nu K_\nu\).
In particular \(r\) increases with \(\nu\). At order \(\tfrac12\),
\(v=\sinh(u/2)\) gives \(K_{1/2}(z)=\sqrt{\pi/(2z)}\e^{-z}\), whence
\(r\le z+1/2\). For \(v\ge z>0\), comparison of the two integrands in~\eqref{eq:B}
gives \(K_\nu(v)\le\e^{-(v-z)}K_\nu(z)\), since \(\cosh u\ge1\).
In particular,
\[
 \int_z^\infty K_\nu(v)^2\,\frac{\dd v}{v}
 \le \frac{K_\nu(z)^2}{z}\int_z^\infty \e^{-2(v-z)}\dd v
 =\frac{K_\nu(z)^2}{2z}.
\]
Inserting this bound into the formula for \(p\) gives
\begin{equation}\label{eq:r}
 z<r\le z+\tfrac12,\qquad 0\le p\le\nu/z.
\end{equation}

Put \(E=(L^2B_\nu)/K_\nu\). In the basis \(K_\nu,DK_\nu\), \(L\) takes
coefficients \((b,c)\) to
\((Db+5b/2+(z^2+\nu^2)c,\ b+Dc+5c/2)\).
Starting from \((z^2+9,6)\) in~\eqref{eq:B}, applying this rule twice yields
\[
 E=z^4+\tfrac{285}{4}z^2+\tfrac{225}{4}
       -(15z^2+\tfrac{165}{2})r+\nu^2(z^2+39-6r).
\]
The coefficient of \(\nu^2\) is at least \(z^2-6z+36>0\). Thus
\[
 4E\ge P(z):=4z^4-60z^3+255z^2-330z+60.
\]
At \(z=25/2\) the values of \(P,P',P'',P'''\) are respectively
\(32495/2,9170,3510,840\), and \(P^{(4)}=96\). This proves \(E>0\)
for \(z\ge4\pi>25/2\). For the order derivative, \eqref{eq:r} gives
\begin{align*}
 \partial_\nu E
 &=2\nu(z^2+39-6r)-(15z^2+165/2+6\nu^2)p\\
 &\ge \nu(2z^2-27z+72-84/z).
\end{align*}
The parenthesis is \(1007/25\) at \(25/2\) and has positive derivative
thereafter. Finally \(K_\nu>0\), and \(\partial_\nu K_\nu>0\) for
\(\nu>0\) by~\eqref{eq:B}. Differentiating \(K_\nu E\) proves the lemma.
\end{proof}

\begin{proposition}\label{pr:shapes}
For \(k\ge2\) and \(0\le\tau\le\tfrac12\), the profile \(H_{\tau,k}\)
is positive, strictly decreasing and strictly convex on \([0,\infty)\);
its Fourier transform is strictly positive on \(\R\). The same assertions
hold for \(Q_{\tau,k}\) when \(0<\tau\le\tfrac12\).
For each fixed \(N\ge1\) and \(0<\tau\le\tfrac12\), \(\cQ_{\tau,N}(t)<0\)
for all sufficiently large \(|t|\).
\end{proposition}
\begin{proof}
For \(k\ge2,\ x\ge0\) we have \(2\pi k\e^x\ge4\pi\).
Lemma~\ref{lem:conv} gives strictly positive second half-line derivatives
for \(G_{\tau,k}\), and also for \(\partial_\tau G_{\tau,k}\) at positive
tilt. The functions and their first derivatives tend to zero at infinity by
the integral estimates of Subsection~\ref{sec:modes}. Integrating back from
infinity gives strict decrease and positivity, with strictly negative right
derivative at zero. In \(a_k(\tau)/k^2\), each pair of distinct divisors
\(j,k/j\) contributes \(2\cosh(\tau\log(k/j^2))\). If \(k\) is a
square, the divisor \(\sqrt k\) contributes \(1\).
Thus \(a_k>0\) and \(a_k'>0\) when \(k\ge2,\ \tau>0\). Formula~\eqref{eq:mode}
gives the assertions about the profiles.

For any one of these even profiles \(G\), two half-line integrations by
parts give the strict version of P\'olya's convexity criterion \cite{Polya}:
\begin{equation}\label{eq:polya}
 \wt G(t)=\frac{2}{t^2}\int_0^\infty G''(x)(1-\cos(tx))\dd x>0
 \qquad(t\ne0).
\end{equation}
Indeed \(\int_0^\infty G''=-G'(0+)\), and \(G''>0\) on the open half-line.
At \(t=0\) positivity follows by integration of \(G\).

The complete \(Q_\tau\) is smooth and even by~\eqref{eq:HQ}; individual
even extensions can have a corner at zero. Write \(S=\sum_{k\le N}Q_{\tau,k}\).
By absolute convergence of the right derivatives,
\(S'(0+)=-\sum_{k>N}Q_{\tau,k}'(0+)>0\).
Four integrations by parts on \([0,\infty)\), using \(S^{(4)}\in L^1\)
there, yield
\[
 \cQ_{\tau,N}(t)=-2S'(0+)t^{-2}+O_{\tau,N}(t^{-4}).
\]
This proves the last assertion, including \(N=1\). Since the mode added at
each successive cutoff has strictly positive transform, the initial interval
of nonnegativity can only expand. Its finiteness and monotonicity follow.
\end{proof}

\section{Size of the omitted modes}\label{sec:tails}

The arithmetic weight and the Fourier profile play different roles. The
Fourier profile has a common limit when frequency and mode index have
comparable size. The divisor weight persists after passage to that limit.

\begin{lemma}\label{lem:scale}
Uniformly in \(0\le\nu\le\tfrac12\) and real \(t\), as \(k\to\infty\),
\begin{equation}\label{eq:scale}
 \e^{2\pi k}k^{-1/2}\wt G_{\nu,k}(t)
   =\frac{8\pi^3}{4\pi^2+(t/k)^2}+o(1).
\end{equation}
The error in~\eqref{eq:scale} is absolute. For every \(k\ge2\), uniformly
in \(\nu,t\),
\begin{equation}\label{eq:L1G}
 |\wt G_{\nu,k}(t)|+|\wt{\partial_\nu G_{\nu,k}}(t)|
                      \ll k^{1/2}\e^{-2\pi k}.
\end{equation}
\end{lemma}
\begin{proof}
From~\eqref{eq:Kbound}, \(B_\nu\) and \(\partial_\nu B_\nu\) are bounded
by \(Cz^{3/2}\e^{-z}\) on \(z\ge1\). The corresponding bound for the two
profiles in the variable \(x\ge0\) is
\(Ck^{3/2}\e^{4x}\e^{-2\pi k\e^x}\).
The inequality \(\e^x\ge1+x\) pays a factor \(O(1/k)\) after integration,
giving~\eqref{eq:L1G}.

For fixed \(v\ge0\), \eqref{eq:Kasym} in~\eqref{eq:B} gives
\[
 g_{\nu,k}(v):=k^{-3/2}\e^{2\pi k}G_{\nu,k}(v/k)
      \longrightarrow 2\pi^2\e^{-2\pi v}.
\]
Here \(k(\e^{v/k}-1)\to v\) and
\((2\pi)^{3/2}\sqrt{\pi/2}=2\pi^2\). The convergence is uniform in \(\nu\)
at fixed \(v\). The preceding profile bound gives
\(|g_{\nu,k}(v)|\le C\e^{4v/k-2\pi v}\le C\e^{-(2\pi-2)v}\).
Dominated convergence applied to the supremum of the difference over \(\nu\)
therefore gives uniform \(L^1(\R_+)\) convergence. In the even Fourier
integral the change \(x=v/k\) leaves the factor \(\cos((t/k)v)\) of
modulus at most one. As a result the error is uniform for all real \(t\).
The limiting cosine integral equals
\[
 4\pi^2\int_0^\infty \e^{-2\pi v}\cos((t/k)v)\dd v
       =\frac{8\pi^3}{4\pi^2+(t/k)^2}.
 \qedhere
\]
\end{proof}

Write \(\sigma=\tfrac12+\tau\) and \(U_{\tau,N}=\sum_{k>N}\wt Q_{\tau,k}\).
For \(\tau>0\) this tail is strictly positive. In particular,
\begin{equation}\label{eq:subtract}
             \cQ_{\tau,N}(t)=W_\tau(t)-U_{\tau,N}(t).
\end{equation}

\begin{lemma}\label{lem:U}
For \(N\ge2\), \(0<\tau\le\tfrac12\) and all real \(t\),
\begin{align}
 U_{\tau,N}(t)&\ll_\epsilon \e^{-2\pi N}N^{\sigma+2+\epsilon}
                                  &&(\epsilon>0),\label{eq:Uup}\\
 U_{1/2,N}(t)&\ll \e^{-2\pi N}N^3\log^2(2N).\label{eq:Uendpoint}
\end{align}
For fixed \(0<a<\tfrac12\), all sufficiently large \(N\), \(a\le\tau\le\tfrac12\),
and \(3N\le t\le5N\), one has
\begin{equation}\label{eq:Ulow}
 U_{\tau,N}(t)\gg_a \e^{-2\pi N}N^{\sigma+2}\log N .
\end{equation}
\end{lemma}
\begin{proof}
From~\eqref{eq:ak},
\(a_k(\tau)\le k^{2+\tau}d(k)\) and \(0\le a_k'(\tau)\le a_k(\tau)\log k\).
For each \(\eta>0\), \(d(k)\ll_\eta k^\eta\): for sufficiently large primes
\(p\), \(j+1\le2^j\le p^{\eta j}\), and the finitely many smaller primes
cost at most \(\prod_{p<2^{1/\eta}}\sup_{j\ge0}(j+1)p^{-\eta j}<\infty\).
Apply~\eqref{eq:L1G} to~\eqref{eq:mode}, choosing \(\eta\) small enough
to absorb \(\log k\) in the allowed \(k^\epsilon\). At the endpoint
\(a_k(1/2)=k^{5/2}\sum_{j\mid k}1/j\le k^{5/2}(1+\log k)\).
Summing the two estimates gives~\eqref{eq:Uup} and~\eqref{eq:Uendpoint}.
Indeed, on writing \(k=N+r\), \(r\ge1\), any fixed polynomial growth in
\(1+r\) is summable against \(\e^{-2\pi r}\), uniformly in the bounded
range of \(\sigma\).

For the lower bound keep \(k=N+1\). The paired divisors \(1,k\) alone,
with the other pairs nonnegative, give
\[
 a_k'(\tau)\ge k^2(k^\tau-k^{-\tau})\log k
              \gg_a k^{2+\tau}\log k .
\]
The transform of \(\partial_\tau G_{\tau,k}\) is positive by~\eqref{eq:polya};
we may keep only \(8\pi^2 a_k'(\tau)\wt G_{\tau,k}\). On the indicated
frequency interval the main term in~\eqref{eq:scale} has a strictly positive
lower bound. Multiplication proves~\eqref{eq:Ulow}, with the absolute
factor \(\e^{-2\pi}\) from the shift to \(N+1\).
\end{proof}

\section{The interval about a critical-line zero}\label{sec:local}

At bounded frequencies the rational factor in~\eqref{eq:scale} approaches
\(2\pi\). This determines the scale of the negative interval for the
untilted transform, including its arithmetic weight.

\begin{lemma}\label{lem:Htail}
For every fixed bounded interval \(J\subset\R\), uniformly in \(t\in J\),
\begin{equation}\label{eq:Htail}
 V_N(t):=\sum_{k>N}\wt H_{0,k}(t)=A_N(1+o_J(1)).
\end{equation}
Moreover
\begin{equation}\label{eq:ANlog}
 \log A_N=-2\pi N+\tfrac52\log N+o(\log N).
\end{equation}
\end{lemma}
\begin{proof}
At zero tilt \(a_k(0)=k^2d(k)\). Lemma~\ref{lem:scale} gives
\(\wt G_{0,k}(t)=2\pi k^{1/2}\e^{-2\pi k}(1+o_J(1))\) uniformly on \(J\).
After multiplication by \(8\pi^2 k^2d(k)\), this is the summand of
\(A_N\) with relative error tending to zero. The supremum of those
errors for \(k>N,\ t\in J\) tends to zero as well. The weights in \(A_N\)
are positive, which proves~\eqref{eq:Htail}. The single term \(k=N+1\)
gives \(A_N\gg N^{5/2}\e^{-2\pi N}\). The bound \(d(k)\ll_\epsilon k^\epsilon\)
gives \(A_N\ll_\epsilon N^{5/2+\epsilon}\e^{-2\pi N}\) for every \(\epsilon>0\).
Taking logarithms gives~\eqref{eq:ANlog}.
\end{proof}

\begin{proof}[Proof of Theorem~\ref{th:zero}]
By~\eqref{eq:completion},
\(\cH_{0,N}(t)=2|\xi(\tfrac12+\ii t)|^2-V_N(t)\).
Every omitted transform is strictly positive by Proposition~\ref{pr:shapes},
so \(\cH_{0,N}(\gamma)<0\) for every \(N\). The multiplicity gives, for real
\(h\) close to zero,
\begin{equation}\label{eq:zerojet}
 2|\xi(s_\gamma+\ii h)|^2=C_\gamma |h|^{2m}
                              \{1+O_\gamma(|h|)\}.
\end{equation}
Choose a fixed bounded interval \(J\) with \(\gamma\) in its interior, and put
\(b_N=(A_N/C_\gamma)^{1/(2m)}\). This tends to zero. For each fixed
\(0<\delta<1\), the complete transform on \(|h|\le(1-\delta)b_N\)
is at most \(A_N(1-\delta)^{2m}(1+o(1))\). By~\eqref{eq:Htail} the
finite transform is strictly negative throughout that interval for large
\(N\). At each of the two points \(h=\pm(1+\delta)b_N\), its value divided
by \(A_N\) tends to \((1+\delta)^{2m}-1>0\).

By continuity the negative component containing \(\gamma\) is bounded
between these two positive points. On either side its extent is at least
\((1-\delta)b_N\) and at most \((1+\delta)b_N\). Letting \(\delta\) decrease
to zero gives both assertions in~\eqref{eq:width}. Now
\(\log |I_N|=(\log A_N-\log C_\gamma)/(2m)+\log2+o(1)\), and~\eqref{eq:ANlog}
gives~\eqref{eq:widthlog}. The use of connected components also allows
other sign changes outside the two positive comparison points.
\end{proof}

\section{Bounds for the complete transform}\label{sec:complete}

There are two estimates to prove. At the endpoint we need two-sided control
up to powers of \(\log t\). For the interior result we need two-sided control
up to \(t^\epsilon\), uniformly on each fixed interval of tilts. We give
the analytic details to specify exactly where RH enters the second estimate.

\subsection{The real logarithmic derivative}

For \(s=\sigma+\ii t\) with \(\xi(s)\ne0\), taking real parts in the classical
partial fraction expansion of \(\xi'/\xi\)~\cite[p.~193, eq.~(3)]{MSZ} gives
\begin{equation}\label{eq:poisson}
 L_\sigma(t):=\re\frac{\xi'}{\xi}(\sigma+\ii t)
     =\sum_{\beta+\ii\gamma}
         \frac{\sigma-\beta}{(\sigma-\beta)^2+(t-\gamma)^2}.
\end{equation}
The sum counts the zeros \(\rho=\beta+\ii\gamma\) with multiplicity and
converges absolutely since \(0\le\beta\le1\) and
\(\sum_\rho(1+\gamma^2)^{-1}<\infty\). For \(\sigma=\tfrac12+\tau\),
\eqref{eq:completion} gives \(W_\tau(t)=4|\xi(s)|^2L_\sigma(t)\).

For \(b>0\) and \(u\in\R\), write \(P_b(u)=b/(b^2+u^2)\).
On the comparison line \(\sigma=2\), every summand in~\eqref{eq:poisson}
is \(P_{2-\beta}(t-\gamma)>0\), so \(L_2(t)>0\) for all real \(t\). Applying Stirling's
formula to the gamma factor and using
\((\zeta'/\zeta)(2+\ii t)=O(1)\) from the Euler product, we obtain
\begin{equation}\label{eq:L2}
 L_2(t)=\tfrac12\log(t/(2\pi))+O(1)\qquad(t\to+\infty).
\end{equation}

\subsection{The endpoint, without RH}

\begin{lemma}\label{lem:endline}
\(W_{1/2}(t)>0\) on the entire real line. With an absolute \(B>0\) and all
sufficiently large \(t\),
\begin{equation}\label{eq:Wend}
 t^4 \e^{-\pi t/2}(\log t)^{-B}\ll W_{1/2}(t)
                         \ll t^4\e^{-\pi t/2}(\log t)^B .
\end{equation}
\end{lemma}
\begin{proof}
For real \(h>0\), nonnegativity of \(3+4\cos v+\cos(2v)=2(1+\cos v)^2\)
in the Euler product gives
\begin{equation}\label{eq:trig}
 1\le\zeta(1+h)^3|\zeta(1+h+\ii t)|^4|\zeta(1+h+2\ii t)|.
\end{equation}
Indeed the logarithm of the right-hand side is a sum of nonnegative terms
over prime powers with \(v=jt\log p\). If \(\zeta(1+\ii t)=0\) at a fixed
\(t\ne0\), the right side tends to zero as \(h\downarrow0\): the pole of
\(\zeta(1+h)\) is simple and \(\zeta(1+h+2\ii t)\) stays bounded.
The simple pole with residue one also gives \(\xi(1)=1/2\).
These facts about the pole follow, for instance, from~\eqref{eq:EM} below.

We have \(L_1(t)\ge L_2(t)/6\). By the functional equation and conjugation,
each zero \(\rho=\beta+\ii\gamma\) with \(\beta<1/2\) is paired with
\(1-\bar\rho=1-\beta+\ii\gamma\), a zero of the same multiplicity.
Set \(x=1-\beta\) and \(u=t-\gamma\). For \(b=2-\beta\) or \(1+\beta\),
\(x\le b\) and \(x/b\ge1/3\). Hence
\(P_x(u)\ge P_b(u)/3\), and \(P_x\) alone dominates a sixth of the sum
of the two line-\(2\) terms. The other line-\(1\) term is nonnegative.
Each critical-line zero satisfies the comparison on its own and is counted
with its multiplicity. Summing~\eqref{eq:poisson} proves the assertion,
including positivity for every \(t\). Equation~\eqref{eq:completion} gives
\(W_{1/2}(t)=4|\xi(1+\ii t)|^2L_1(t)>0\).

For quantitative bounds, let \(X\) be a positive integer. On \(\re s>0\),
\(s\ne1\),
\begin{equation}\label{eq:EM}
 \zeta(s)=\sum_{n\le X}n^{-s}+\frac{X^{1-s}}{s-1}
                  -s\int_X^\infty\{u\}u^{-s-1}\dd u .
\end{equation}
For \(\re s>1\) this follows by integrating the counting function
\(\lfloor u\rfloor-X\) for the tail and splitting off \(u-X\);
absolute convergence of the integral gives the stated continuation.
For heights \(t,2t\) with \(t\) large, take \(X=\lfloor t^2\rfloor\) and
\(\ell=\log t\). In \(1\le\re s\le2\) the finite sums for \(\zeta,\zeta'\)
are \(O(\ell),O(\ell^2)\), by comparison with \(\sum_{n\le X}1/n\) and
\(\sum_{n\le X}(\log n)/n\). The integral remainder is
\(O(|s|X^{-\re s})\), and its derivative is
\(O((1+|s|)(1+\log X)X^{-\re s})\) on this strip. The explicit pole
term and its derivative are bounded there at these heights. Thus
\(|\zeta|\ll\ell,\ |\zeta'|\ll\ell^2\) uniformly.

Take \(h=\ell^{-16}\) in~\eqref{eq:trig}. Since
\(\zeta(1+h)\ll1/h\), it follows that
\(|\zeta(1+h+\ii t)|\gg h^{3/4}\ell^{-1/4}=\ell^{-49/4}\).
Moving horizontally to \(1+\ii t\) costs \(O(h\ell^2)=O(\ell^{-14})\).
In particular \(|\zeta(1+\ii t)|\gg\ell^{-13}\).
Thus \(\lvert(\zeta'/\zeta)(1+\ii t)\rvert\ll\ell^{15}\).
Logarithmic differentiation of the defining formula for \(\xi\) and
Stirling's formula give \(L_1(t)\ll\ell^{15}\). The comparison
\(L_1(t)\ge L_2(t)/6\) and~\eqref{eq:L2} give \(L_1(t)\gg\ell\).
Finally, uniformly for \(1/2\le\sigma\le1,\ t\) large, Stirling gives
\begin{equation}\label{eq:Stir}
 |\xi(\sigma+\ii t)|^2
       \asymp t^{\sigma+3}\e^{-\pi t/2}|\zeta(\sigma+\ii t)|^2 .
\end{equation}
The power comes from \(t^4\) in \(|s(s-1)|^2\) and \(t^{\sigma-1}\) in
\(|\Gamma(s/2)|^2\). Using \(\sigma=1\) in this formula proves
\eqref{eq:Wend}, for example with \(B=30\).
\end{proof}

\subsection{A compact interval of positive tilts under RH}

\begin{lemma}\label{lem:strip}
Assume RH and fix \(0<a<\tfrac12\). Then \(W_\tau(t)>0\) for all real \(t\)
and \(a\le\tau\le\tfrac12\). For every \(\epsilon>0\), uniformly in these
tilts at all sufficiently large \(t\), with \(\sigma=\tfrac12+\tau\),
\begin{equation}\label{eq:Wstrip}
 t^{\sigma+3-\epsilon}\e^{-\pi t/2}
       \ll_{a,\epsilon} W_\tau(t)
       \ll_{a,\epsilon} t^{\sigma+3+\epsilon}\e^{-\pi t/2}.
\end{equation}
\end{lemma}
\begin{proof}
Under RH, \(\beta=\tfrac12\) and \(\sigma-\beta=\tau\), so each summand
in~\eqref{eq:poisson} is
\mbox{\(P_\tau(t-\gamma)=\frac{\tau}{\tau^2+(t-\gamma)^2}>0\)}.
For \(0<v\le b\), \(v/b\le P_v(u)/P_b(u)\le b/v\).
Taking \(v=\tau\) and \(b=3/2\) and summing over the zeros gives
\begin{equation}\label{eq:Pcomp}
       (2a/3)L_2(t)\le L_{\sigma}(t)\le (3/(2a))L_2(t).
\end{equation}
Thus \(L_\sigma(t)>0\) everywhere and \(L_\sigma(t)\asymp_a\log t\) at large
heights. The zero-free half-plane also gives \(\xi(\sigma+\ii t)\ne0\)
on these lines. This proves the positivity assertion.

For clarity we give a proof of the subpower bounds for \(\zeta\) and its
reciprocal in this range, a classical consequence of RH
\cite[\S14.2]{Titchmarsh}. At a large fixed \(t\), use the disk centered at
\(z_t=2+\ii t\) of radius \(R=3/2-a/4\). The disk and a neighborhood of
its closure are free of zeros and of the pole. Choose an analytic logarithm
\(g_t=\log\zeta\) there, agreeing at the center with the absolutely convergent
Euler logarithm. In \(|z-z_t|\le1/2\) that series gives \(|g_t|\ll1\).
On the entire disk, \eqref{eq:EM} with \(X=\lfloor t^2\rfloor\) gives
\(|\zeta(z)|\ll_a (2+t)^C\) for an absolute \(C\), since
\(\re z\ge1/2+a/4\). Thus \(\re g_t\le C_a\log t\) there.
On the smaller disk of radius \(R_1=3/2-a/2\), it follows that
\(\max|g_t|\ll_a\log t\). To deduce this from the bound on \(\re g_t\),
write \(h=g_t-g_t(z_t)\), \(M=1+\max_{|z-z_t|\le R}\re h\). The function
\(h/(2M-h)\) maps this disk into the unit disk and vanishes at its center.
Schwarz's lemma bounds its modulus by \(R_1/R\) on the smaller disk,
which gives the asserted bound for \(g_t\).

Put \(r=3/2-a\). The three-circles inequality, or the maximum principle for
\(\log|g_t|\) on \(1/2\le|z-z_t|\le R_1\) with a majorant affine in
\(\log|z-z_t|\), gives
\[
 \max_{|z-z_t|\le r}|g_t(z)|\ll_a(\log t)^{\theta_a},\qquad
 \theta_a=\frac{\log(r/(1/2))}{\log(R_1/(1/2))}<1 .
\]
The intended segment \(1/2+a\le\sigma\le1\) lies in this disk at height
\(t\). Thus both \(|\zeta(\sigma+\ii t)|\) and its reciprocal are at most
\(\exp(C_a(\log t)^{\theta_a})\). Substitution in~\eqref{eq:Stir}, followed
by~\eqref{eq:Pcomp} and \(W_\tau=4|\xi|^2L_\sigma\), proves~\eqref{eq:Wstrip}:
the square of this subpower factor and the logarithm are absorbed by
\(t^\epsilon\) for each prescribed \(\epsilon>0\).
\end{proof}

\section{Locating the first sign loss}\label{sec:window}

\begin{proof}[Proof of Theorem~\ref{th:window} and \eqref{eq:inverse}]
Finiteness and monotonicity follow from Proposition~\ref{pr:shapes}.
We use~\eqref{eq:subtract} for the quantitative assertions.

\smallskip\noindent\textit{Endpoint.}
With \(B\) from Lemma~\ref{lem:endline}, choose \(D\) so that
\(\pi D/2>B+3\), and put
\[
 t_\pm=4N+(2/\pi)\log N\pm D\log\log N .
\]
Both \(t_\pm/N\to4\), and
\[
 \e^{-\pi t_\pm/2}
      =\e^{-2\pi N}N^{-1}(\log N)^{\mp\pi D/2}.
\]
At \(t_+\), the upper bound in~\eqref{eq:Wend} and the lower bound
\eqref{eq:Ulow} imply
\[
 W_{1/2}(t_+)/U_{1/2,N}(t_+)
       \ll(\log N)^{B-\pi D/2-1}\longrightarrow0.
\]
For the opposite comparison, fix \(T_0\) large enough for~\eqref{eq:Wend}
and for \(t^4\e^{-\pi t/2}(\log t)^{-B}\) to decrease on
\([T_0,\infty)\). Its logarithmic derivative is
\(4/t-\pi/2-B/(t\log t)\). On the entire interval \(T_0\le t\le t_-\),
\[
 W_{1/2}(t)\gg \e^{-2\pi N}N^3(\log N)^{\pi D/2-B}.
\]
This exceeds the uniform upper bound~\eqref{eq:Uendpoint} for large \(N\).
On \([0,T_0]\), \(W_{1/2}\) has a strictly positive minimum and the tail
tends uniformly to zero. Thus \(\cQ_{1/2,N}>0\) on \([0,t_-]\), while
\(\cQ_{1/2,N}(t_+)<0\). The resulting bracket for \(\rho_N(1/2)\)
proves~\eqref{eq:endpoint}.

\smallskip\noindent\textit{Uniformity in the tilt.}
Assume RH. Fix \(a\) as in the theorem and an arbitrary \(\delta>0\).
Take \(0<\epsilon<\min(1,\pi\delta/4)\) and put
\[
 u_\pm=4N+\bigl(2/\pi\pm\delta\bigr)\log N .
\]
Then \(\e^{-\pi u_\pm/2}=\e^{-2\pi N}N^{-1\mp\pi\delta/2}\).
For all \(a\le\tau\le\tfrac12\), the quotient at the upper point satisfies
\[
 \frac{W_\tau(u_+)}{U_{\tau,N}(u_+)}
       \ll_{a,\epsilon,\delta}
       \frac{N^{\epsilon-\pi\delta/2}}{\log N}\longrightarrow0
\]
by~\eqref{eq:Wstrip} and~\eqref{eq:Ulow}. For \(t\ge T_0(a,\epsilon)\)
the function \(t^{\sigma+3-\epsilon}\e^{-\pi t/2}\) decreases
simultaneously for \(1/2+a\le\sigma\le1\), since
\((\sigma+3-\epsilon)/t-\pi/2<0\). Hence, throughout \(T_0\le t\le u_-\),
\[
 W_\tau(t)\gg_{a,\epsilon,\delta}
   \e^{-2\pi N}N^{\sigma+2+\pi\delta/2-\epsilon}.
\]
This dominates~\eqref{eq:Uup} uniformly, since \(\pi\delta/2>2\epsilon\).
Joint continuity and strict positivity of \(W_\tau\) give a positive minimum
on \([a,1/2]\times[0,T_0]\), and~\eqref{eq:Uup} tends to zero uniformly
there. We have proved \(u_-\le\rho_N(\tau)\le u_+\) for every tilt in the
interval and all sufficiently large \(N\). Since \(\delta\) was arbitrary,
this is~\eqref{eq:fixedstrip}.

\smallskip\noindent\textit{Inversion.}
Let \(x(T)=T/4-(\log T)/(2\pi)\). For
\(n_\pm=\lfloor x(T)\pm M\log\log T\rfloor\), one has
\[
 4n_\pm+(2/\pi)\log n_\pm
        =T\pm4M\log\log T+O(1).
\]
Choose \(M\) larger than the error constant in~\eqref{eq:endpoint}.
Then \(\rho_{n_-}(1/2)<T<\rho_{n_+}(1/2)\) for large \(T\).
The positivity interval includes its finite right endpoint by continuity
whenever it has positive length. Monotonicity in \(N\) therefore gives
\(n_-<N_{\min}(T)\le n_+\), proving~\eqref{eq:inverse}.
\end{proof}

At \(t=4N\) the exponentials \(\e^{-2\pi N}\) and \(\e^{-\pi t/2}\)
are equal. For \(3N\le t\le5N\) the complete-transform bounds
carry the power \(N^{\sigma+3}\), and the arithmetic tail carries
\(N^{\sigma+2}\), up to subpower factors. Raising \(t\) by
\((2/\pi)\log N\) compensates for this one power of \(N\).
At a fixed critical-line zero the same exponential \(\e^{-2\pi N}\)
is balanced instead by the local power \(|t-\gamma|^{2m}\). This accounts
for the different scales in Theorems~\ref{th:window} and~\ref{th:zero}.

\section*{Scope of the formal verification}\label{sec:formal}
\addcontentsline{toc}{section}{Scope of the formal verification}
Drawing on their joint research, Wei Sun and OpenAI's ChatGPT (Alex,
GPT-6-based Codex) prepared and revised this manuscript from 25 through
29 November 2026. Sun's prompts asked for self-contained proofs of what the finite arithmetic
sums reveal about \(\xi\), even when RH is assumed. He also asked Alex to ``complete Lean
verification'' of the results in this paper. Alex wrote the proof scripts,
which the Lean kernel checked. Symbolic algebra was used to check the polynomial
coefficients and derivatives in Lemma~\ref{lem:conv}.

All nine numbered results in this paper have formal proofs.
Lean 4.24.0 proves Theorem~\ref{th:zero}, Proposition~\ref{pr:shapes} and
Lemmas~\ref{lem:conv}, \ref{lem:scale}, \ref{lem:U} and~\ref{lem:Htail} for the
functions defined in this paper. At a critical-line zero it proves
strict negativity for every positive integer cutoff, the two endpoint
equivalents~\eqref{eq:width}, and the logarithmic length
formula~\eqref{eq:widthlog}. The only zero hypothesis specifies the
chosen zero and its multiplicity. The formalization also proves the
Mellin and correlation formulas, including~\eqref{eq:completion} at zeros.

For Lemma~\ref{lem:conv} the formal proof uses a comparison at two
Bessel orders. With \(r_\nu=-zK_\nu'/K_\nu\), the Bessel equation gives
\[
 (r_\beta-r_\alpha) K_\alpha(z)K_\beta(z)
 = (\beta^2-\alpha^2)\int_z^\infty K_\alpha(v)K_\beta(v)\,\frac{\dd v}{v}.
\]
Indeed, the derivative of
\((vK'_\alpha)K_\beta-K_\alpha(vK'_\beta)\) equals
\((\alpha^2-\beta^2)K_\alpha K_\beta/v\), and the boundary term at
infinity vanishes by exponential decay. The integral is at most
\(K_\alpha(z)K_\beta(z)/(2z)\). Difference quotients therefore give
\(0\le\partial_\nu r_\nu\le\nu/z\) for \(\nu>0\).
Together with the half-order value this proves the two strict curvature
inequalities on the full stated domain.

Lean also proves positivity, strict decrease and strict convexity on
the right half-line and strict Fourier positivity for the original
\(H_{\tau,k}\) and \(Q_{\tau,k}\) modes throughout their respective tilt
ranges in Proposition~\ref{pr:shapes}. For the order derivative it takes
one-sided limits of differences of Bessel profiles, including at
\(\tau=1/2\). For Lemma~\ref{lem:U}, one divisor pair in the first
omitted mode supplies the signed lower bound throughout \(3N\le t\le5N\).
Summing the absolute bounds gives the upper estimates, including the
endpoint factor \(\log^2(2N)\).

For the fixed-cutoff sign in Proposition~\ref{pr:shapes}, fix
\(k=N+1\). In~\eqref{eq:polya}, \(G_{\tau,k}''\) has a positive
minimum on \([1,2]\), and
\(\int_1^2(1-\cos(tx))\dd x\ge1-2/|t|\).
Since \(a_k'(\tau)>0\), the first omitted mode contributes at least
\(c_{N,\tau}|t|^{-2}\) when \(|t|\ge4\); every later mode is positive.
The Mellin integral bounds \(\xi\) on \(|\re s|\le3\).
On \(\re s=2\), six gamma recurrences, the inequality
\(|\Gamma(z+6)|\le\Gamma(\re z+6)\) for \(\re z>0\), and the
absolutely convergent zeta series bound \([s(1-s)]^2\xi(s)\).
The functional equation gives the same bound on \(\re s=-1\).
Polynomial growth and the Phragm\'en--Lindel\"of theorem extend the
bound throughout the strip. Cauchy's estimate on unit circles in
\(|\re s|\le3\) also gives \(\xi'(s)=O(1)\) for \(|\re s|\le2\).
Thus \(W_\tau(t)=O(|t|^{-4})\) by~\eqref{eq:completion}.
Subtraction gives the finite-sum sign for every \(N\ge1\).
The finiteness and monotonicity of \(\rho_N\) are also verified.

Lean 4.31.0 proves Lemma~\ref{lem:endline} for \(W_{1/2}\), with \(B=14\).
The formal proof uses the canonical product~\cite{Harr}, counting zeros
with multiplicity, to obtain the absolutely convergent real sum
in~\eqref{eq:poisson} on the lines \(\sigma=1,2\). Pairing zeros in
that sum gives \(L_1(t)\ge L_2(t)/6>0\) at every real height.
The digamma series and the absolutely convergent series for
\(\zeta'/\zeta\) give \(L_2(t)=\tfrac12\log t+O(1)\).
Euler's reflection formula gives, for \(t\ne0\),
\[
 |\xi(1+\ii t)|^2
  =\frac{(1+t^2)t^2|\zeta(1+\ii t)|^2}{4\cosh(\pi t/2)} .
\]
The estimates of \cite{PNTlean} bound \(|\zeta|\), \(|\zeta|^{-1}\) and
\(|\zeta'/\zeta|\) on this line by constant multiples of \(\log t\),
\((\log t)^7\) and \((\log t)^9\), respectively, for \(t>3\).
Multiplication gives the two bounds in~\eqref{eq:Wend}.

Lean 4.31.0 also proves Lemma~\ref{lem:strip} under RH. Write \(z_t=2+\ii t\)
and let \(g(w)\) be the analytic logarithm of \(\zeta(z_t+w)/\zeta(z_t)\)
on a disk about \(w=0\), normalized by \(g(0)=0\).
Borel--Carath\'eodory bounds \(g\) by a constant on an inner disk
and by \(O_a(\log t)\) on a larger one. On their intervening annulus
the maximum principle applied to \(g(w)/w^n\), for fixed large \(n\),
makes the outer contribution smaller than any prescribed multiple of
\(\log t\) on the target strip. This proves the uniform subpower bounds
for both \(\zeta\) and its reciprocal. The horizontal derivative of
\(\log|\Gamma((\sigma+\ii t)/2)|^2\) is the real part of \(\Gamma'/\Gamma\)
at \((\sigma+\ii t)/2\). The mean value theorem then gives the gamma estimate
by comparison with \(\sigma=1\).
Combining these estimates with the Poisson comparison proves both
inequalities in~\eqref{eq:Wstrip} and positivity for every real height.

For Theorem~\ref{th:window}, Lean 4.31.0 proves the identity
\(\cQ_{\tau,N}=W_\tau-U_{\tau,N}\) directly for the original finite
arithmetic sum. On \(0\le t\le3N\), a positive minimum on bounded
heights and the large-height estimates give
\(W_\tau(t)\ge\e^{-3\pi N/2}\) for sufficiently large \(N\).
The global upper bound for the omitted tail is smaller. On
\(3N\le t\le5N\), subtracting the logarithms of the two positive
quantities gives the leading term
\[
 \frac{\pi}{2}\left(4N+\frac{2}{\pi}\log N-t\right).
\]
The error is \(O(\log\log N)\) at the endpoint and \(o(\log N)\)
uniformly on each stated positive-tilt strip under RH. Positivity through
the lower bracket and strict negativity at the upper test point give
\eqref{eq:endpoint} and~\eqref{eq:fixedstrip}. Two integer bounds on the least
cutoff yield the unconditional inverse formula~\eqref{eq:inverse}.

\section*{Contributions and use of generative AI}
\addcontentsline{toc}{section}{Contributions and use of generative AI}
\begin{samepage}
This research grew out of Wei Sun's work on modulus representations of
\(\xi\) and inequalities equivalent to RH \cite{SunInequality,SunModulus}.
Sun formulated the questions and directed the joint research that he and
OpenAI's ChatGPT (Alex) carried out from 17 August through 25 November
2026. Alex developed
the Bessel estimates, compared the complete transform with its omitted
modes and worked out the local argument at a zero of given multiplicity.
Their joint work also produced Lean proofs of all nine numbered results
(\hyperref[sec:formal]{see the preceding section}). Alex checked proof
steps and drafted the exposition.
Wei Sun reviewed the manuscript, revised it jointly with Alex, and takes
full responsibility for the mathematics and the text.

\end{samepage}


\begingroup
\small
\begin{thebibliography}{99}
\setlength{\itemsep}{2pt plus 0.5pt minus 0.5pt}
\addcontentsline{toc}{section}{References}
\bibitem{Csordas}
G.~Csordas, Fourier transforms of positive definite kernels and the Riemann
\(\xi\)-function, \textit{Comput. Methods Funct. Theory} \textbf{15}
(2015), 373--391. \href{https://doi.org/10.1007/s40315-014-0105-8}
{doi:10.1007/s40315-014-0105-8}.

\bibitem{Haglund}
J.~Haglund, Some conjectures on the zeros of approximates to the Riemann
\(\Xi\)-function and incomplete gamma functions,
\textit{Cent. Eur. J. Math.} \textbf{9} (2011), 302--318.
\href{https://doi.org/10.2478/s11533-010-0095-3}{doi:10.2478/s11533-010-0095-3}.

\bibitem{Harr}
T.~Harr, \textit{Strict Logarithmic-Derivative Criteria for Zero Location
in Lean: Polynomials, Genus-One Products, the Canonical Xi Divisor,
and a Fully Formalized Xi--Lagarias Equivalence}, Lean software,
version~0.1.0. \href{https://github.com/alejandrozu/RiemannHypothesis-Formalization/blob/main/CITATION.cff}{Software citation record}.
Our proof build uses
\href{https://github.com/alejandrozu/RiemannHypothesis-Formalization}{the \texttt{alejandrozu} repository}
at revision
\href{https://github.com/alejandrozu/RiemannHypothesis-Formalization/tree/753938937e624d1c8bbe4c208517660a672d1b8e}{\texttt{753938937e62}}.

\bibitem{MSZ}
Yu.~Matiyasevich, F.~Saidak and P.~Zvengrowski, Horizontal monotonicity
of the modulus of the zeta function, \(L\)-functions, and related functions,
\textit{Acta Arith.} \textbf{166} (2014), 189--200.
\href{https://doi.org/10.4064/aa166-2-4}{doi:10.4064/aa166-2-4}.
Equation~(3) appears without a number in the proof of Theorem~1.1
on p.~4 of \href{https://arxiv.org/abs/1205.2773v1}{arXiv:1205.2773v1}.

\bibitem{DLMF}
F.~W.~J.~Olver et al. (eds.), \textit{NIST Digital Library of Mathematical
Functions}, National Institute of Standards and Technology, Chapter~10,
\href{https://dlmf.nist.gov/10.32}{\S10.32} and
\href{https://dlmf.nist.gov/10.40}{\S10.40}.

\bibitem{Polya}
G.~P\'olya, Remarks on characteristic functions, in \textit{Proceedings
of the Berkeley Symposium on Mathematical Statistics and Probability},
University of California Press, Berkeley, 1949, 115--123.
\href{https://digicoll.lib.berkeley.edu/record/112797/files/math_s1_article-08.pdf}{Berkeley Digital Collections}.

\bibitem{PNTlean}
PrimeNumberTheoremAnd contributors, \textit{PrimeNumberTheoremAnd},
Lean software. Module \texttt{ZetaBounds} at revision
\href{https://github.com/AlexKontorovich/PrimeNumberTheoremAnd/tree/f54d1f4edb45b91da4392df975451427d10c4122}{\texttt{f54d1f4e}}.
Modules \texttt{ZetaBounds} and R.~Sneiderman's
\texttt{DigammaSeries} at revision
\href{https://github.com/AlexKontorovich/PrimeNumberTheoremAnd/tree/b01237e8389be17d52c32419cc8ae8a6d4e4218b}{\texttt{b01237e8}}.

\bibitem{SD}
J.~Sondow and C.~Dumitrescu, A monotonicity property of Riemann's xi
function and a reformulation of the Riemann hypothesis,
\textit{Period. Math. Hungar.} \textbf{60} (2010), 37--40.
\href{https://doi.org/10.1007/s10998-010-1037-3}{doi:10.1007/s10998-010-1037-3}.

\bibitem{SunInequality}
W.~Sun, An equivalent inequality for the Riemann hypothesis, 2023,
revised 2024. \href{https://arxiv.org/abs/2312.12760v5}{arXiv:2312.12760v5}.

\bibitem{SunModulus}
W.~Sun, Modulus representation of the Riemann \(\xi\) function and
polynomial inequalities equivalent to the Riemann hypothesis, 2024.
\href{https://arxiv.org/abs/2404.06642v3}{arXiv:2404.06642v3}.

\bibitem{Titchmarsh}
E.~C.~Titchmarsh, \textit{The Theory of the Riemann Zeta-Function},
2nd ed., revised by D.~R.~Heath-Brown, Clarendon Press, Oxford, 1986.

\bibitem{Zagier}
D.~Zagier, The Rankin--Selberg method for automorphic functions which are
not of rapid decay, \textit{J. Fac. Sci. Univ. Tokyo Sect. IA Math.}
\textbf{28} (1982), 415--437.
\href{https://doi.org/10.15083/00039589}{University of Tokyo repository}.
\end{thebibliography}
\endgroup
\end{document}